\documentclass{article}
\usepackage{amssymb,amsmath}

\newtheorem{theorem}{Theorem}
\newtheorem{corollary}{Corollary}
\newtheorem{lemma}{Lemma}

\newtheorem{conjecture}{Conjecture}

\begin{document}

\title{Congruences for $9$-regular partitions modulo 3}
\author{William J. Keith}

\maketitle

\begin{abstract}

It is proved that the number of 9-regular partitions of $n$ is divisible by 3 when $n$ is congruent to 3 mod 4, and by 6 when $n$ is congruent to 13 mod 16.  An infinite family of congruences mod 3 holds in other progressions modulo powers of 4 and 5.  A collection of conjectures includes two congruences modulo higher powers of 2 and a large family of "congruences with exceptions" for these and other regular partitions mod 3.

\end{abstract}

\section{Introduction}

The $\ell$-regular partitions may be described as those in which parts that are of sizes which are multiples of $\ell$ do not appear.  Their generating function is $$B_\ell(q) = \sum_{n=0}^{\infty} b_\ell(n) q^n = \frac{(q^\ell;q^\ell)_{\infty}}{(q;q)_{\infty}}.$$

In recent years there has been a spate of work (\cite{AHS}, \cite{CuiGu}, \cite{FPenn}, \cite{Sell1}, \cite{Webb} and more) on the divisibility properties of $b_\ell(n)$ in arithmetic progression for various $\ell$, especially on their parity.  Many $B_\ell(q)$ possess numerous progressions of values divisible by 2 and 3, in stark contrast to the usual partition function, where such sequences are unknown and generally held to be almost certainly impossible.  In another contrast, theorems of the form $b_\ell(An+B) \equiv 0$ mod $M$ often have $A$ and $M$ coprime.

Furcy and Penniston in \cite{FPenn} were able to prove infinite families of congruences modulo 3 for many $\ell \equiv 1$ mod 3.  Cui and Gu in \cite{CuiGu} find parity results for $\ell = 2^j$ up to 16. Hirschhorn and Sellers in \cite{Sell1} produce results for $\ell = 5$ as well as 13.  The modulus $\ell = 9$ was studied for its parity by Xia and Yao (\cite{XiaYao}); here we study it modulo 3.

The theorems of this paper are

\begin{theorem} For $k \in \mathbb{N} \bigcup \{0\}$, $b_9(4k+3) \equiv 0$ mod 3
\end{theorem}

\noindent and the related theorem

\begin{theorem} $b_9(16k+13) \equiv 0$ mod 6,\end{theorem}

\noindent the latter being slightly stronger than the $a=2$ case of the general theorem

\begin{theorem} $b_9\left(4^an+\frac{10 \cdot 4^{a-1}-1}{3} \right) \equiv 0$ mod 3.\end{theorem}

One may further conjecture, for higher powers of 2, that

\begin{conjecture} $b_9(32k+13) \equiv 0 $ mod 12
\end{conjecture}

\noindent and 

\begin{conjecture} $b_9(64k+13) \equiv 0 $ mod 24.\end{conjecture}

\noindent \textbf{Remark:} The natural extension to $128n+13$ or $128n+77$ fails.

The original preprint version of this paper conjectured in detail the infinite family of congruences

\begin{theorem} $b_9(5^{2a}n+\frac{5^{2a-2}-1}{3} + 5^{2a-2}k) \equiv 0$ mod 3 for any positive integer $a$ and $k=3$,13,18, or 23.\end{theorem}

This was based on numerical evidence for 

\begin{lemma}\label{NineFiveDissect} $\sum_{n=0}^\infty b_9(5n+3) q^{n} \equiv q \frac{(q^{45};q^{45})_\infty}{(q^5;q^5)_\infty}$ mod 3.\end{lemma}

Shortly after the preprint announcement for this paper, the author was contacted by Ernest Xia and Olivia Yao of Jiangsu University, who were able to verify the conjecture by constructing a complete dissection of $B_9(q)$ which includes the term above.  They intend to publish soon.

In the next section we define notations and give the necessary background on modular forms.  In Section 3 we go through the proofs of Theorems 1, 2, and 3.  In the last section we discuss the conjectures and Theorem 4, and make note of the family of related conjectural congruences for other $\ell$.

\section{Notations and modular form background}

A reader familiar with the notation used in the Introduction, and with the use of eta products, Sturm's theorem and Hecke operators in the partition theoretic literature, may safely skip this section.

We employ the usual notation $$(q^a;q^b)_\infty := \prod_{n=0}^{\infty} (1-q^{a+bn}).$$

When we state that two power series $f(q) = \sum a(n)q^n$ and $g(q) = \sum b(n)q^n$ are congruent mod $M$, we mean that $a(n) \equiv b(n) $ mod $M$ for all $n$.

Modular forms are a class of functions which are integral to many congruence proofs in partition theory; it is sufficient for our purposes to use several theorems regarding their properties without defining them explicitly.  We will refer to the (complex vector space of) modular forms of a given weight $k$, level $N$, and character $\chi$ as $M_k(\Gamma_0(N),\chi)$.

Given $q := e^{2\pi i z}$, Dedekind's eta function is $$\eta(z) = q^{1/24} \prod_{n=1}^\infty (1-q^n).$$

Many generating functions of interest in partition theory may be written as products of eta functions.  A theorem of Gordon, Hughes and Newman (\cite{GordonHughes} and \cite{Newman}) gives sufficient conditions for an eta product to be a modular form.  

\begin{theorem} Let $f(z) = \prod_{\delta \vert N} \eta^{r_\delta} (\delta z)$ with $r_\delta \in \mathbb{Z}$.  If $$\sum_{\delta \vert N} \delta r_\delta \equiv 0 \, \text{ mod } 24 \quad \text{ and } \quad \sum_{\delta \vert N} \frac{N}{\delta} r_\delta \equiv 0 \, \text{ mod } 24,$$ then $f(z)$ is a modular form of weight $k=\frac{1}{2} \sum r_\delta$, level $N$, and character $\chi(d) = \left( \frac{(-1)^k s}{d} \right)$, where $\left(\frac{\cdot}{d}\right)$ is the Jacobi symbol and $s = \prod \delta^{r_\delta}$.
\end{theorem}

One may observe that if the first condition is satisfied but all $\delta$ divide some $N$ for which the second condition does not hold, we may simply multiply $N$ by a sufficient quantity (at worst, 24).  This has only the effect of increasing the work necessary for the next theorem, due to Sturm (\cite{Sturm}), which tells us when the difference of two modular forms is divisible by a given prime:

\begin{theorem} Let $f(z) = \sum_{n=0}^\infty a(n) q^n$ and $g(z) = \sum_{n=0}^\infty b(n) q^n$ both be modular forms in $M_k(\Gamma_0(N),\chi)$ for some positive integer $N$ and character $\chi$.  If $$p \vert (a(n) - b(n)) \quad \forall \quad 0 \leq n \leq \frac{k}{12} N \cdot \prod_{d \, \text{\emph{prime}}: d \vert N} \left(1+\frac{1}{d}\right),$$ then $f(z) \equiv g(z)$ mod $p$.
\end{theorem}

In particular, since 0 is a modular form in every $M_k(\Gamma_0(N),\chi)$, if all coefficients of a modular form $f(z)$ up to the Sturm bound are divisible by $p$, then all coefficients of $f(z)$ are so divisible.

In order to select particular arithmetic progressions to examine with Sturm's theorem it is useful to employ Hecke operators.  For a prime $p$ the Hecke operator $T_p$ (actually $T_{p,k,\chi}$, but the latter two are usually suppressed when these values are clear from context) is an endomorphism on $M_k(\Gamma_0(N),\chi)$.  We write, for $f(z) = \sum_{n=0}^\infty a(n) q^n \in M_k(\Gamma_0(N),\chi)$, $$f \vert T_p := \sum (a(pn) + \chi(p) p^{k-1}a(n/p)) q^n$$ with the convention that $a(n/p) = 0$ for $p \nmid n$.  When $k>1$, the latter term vanishes mod $p$ and we have the factorization property that $$\left( f \cdot \sum_{n=0}^\infty b(n) q^{pn} \right) \vert T_p \equiv \left( \sum_{n=0}^{\infty} a(pn)q^n \right) \left(\sum_{n=0}^\infty b(n) q^n \right) \, \text{ mod } p.$$

\section{Proofs}

\subsection{Proofs of Theorem 1 and 3}

Theorems 1 and 2 rely on the dissection

\begin{lemma}\label{MainLemma} $$B_9(q) = \frac{(q^{12};q^{12})_\infty^2}{(q^{2};q^{2})_\infty^2 {(q^{6};q^{36})_\infty} {(q^{30};q^{36})_\infty}} + q \frac{(q^{12};q^{24})_\infty^2{(q^{36};q^{36})_\infty}}{(q^{4};q^{4})_\infty(q^{4};q^{8})_\infty^6} + 3 q^3 \frac{(q^{24};q^{24})_\infty^2{(q^{36};q^{36})_\infty}}{(q^{4};q^{4})_\infty^3(q^{4};q^{8})_\infty^2}.$$
\end{lemma}

Once this has been conjectured by experimental calculation, the verification is straightforward by means of standard theorems on eta products and congruences for modular forms.  Theorem 2 also requires the classical fact that 

\begin{lemma}\label{InsideThree} $\frac{(q^{12};q^{24})_\infty^2}{(q^4;q^8)_\infty^6} \equiv 1$ mod 3
\end{lemma}

\noindent and showing, via the use of Hecke operators, that

\begin{lemma}\label{EvenLemma} $b_9(16k+13)$ is even. \end{lemma}

We rewrite the identity of Lemma \ref{MainLemma} in terms of eta products, using the identities 

$$\frac{1}{(q,q^{5};q^{6})_\infty} = \frac{(q^{2};q^{2})_\infty (q^{3};q^{3})_\infty}{(q;q)_\infty (q^{6};q^{6})_\infty} \, \text{ and } \, \frac{1}{(q;q^2)_\infty} = \frac{(q^2;q^2)_\infty}{(q;q)_\infty}.$$

We obtain that we wish to prove

\begin{multline}\label{EtaRewrite} q^{-1/3} \frac{\eta(9z)}{\eta(z)} =  q^{-1/3} \frac{\eta(12z)^2 \eta(18z) \eta(12z)}{\eta(2z)^2 \eta(6z) \eta(36z)} \\
+ q^{-1/3} \frac{\eta(8z)^6 \eta(36z) \eta(12z)^2}{\eta(4z)^7 \eta(24z)^2} + 3 q^{-1/3} \frac{\eta(8z)^2\eta(36z)\eta(24z)^2}{\eta(4z)^5}.\end{multline}

Multiplying through by a factor of $q^{1/3} \eta(4z)^4$, we obtain

\begin{multline}\label{EtaRewriteTwo} \frac{\eta(9z) \eta(4z)^4}{\eta(z)} =  \frac{\eta(12z)^2 \eta(18z) \eta(12z)\eta(4z)^4}{\eta(2z)^2 \eta(6z) \eta(36z)} \\
+ \frac{\eta(8z)^6 \eta(36z) \eta(12z)^2}{\eta(4z)^3 \eta(24z)^2} + 3 \frac{\eta(8z)^2\eta(36z)\eta(24z)^2}{\eta(4z)}.\end{multline}

Applying the theorem of Gordon, Hughes and Newman, taking $N=216$, we learn that each term of (\ref{EtaRewriteTwo}) is a modular form of weight 2 on $\Gamma_0 (216)$.  Since the sum of two modular forms of weight $k$ on a given congruence subgroup is again a modular form of weight $k$ for the same congruence subgroup, the right-hand side of (\ref{EtaRewriteTwo}) is a modular form of weight 2 on $\Gamma_0 (216)$ and we may use Sturm's theorem to compare the two sides of the equation.  The bound on coefficients required is $$ \frac{2}{12} 216 \left(1+\frac{1}{2}\right)\left(1+\frac{1}{3}\right) = 72.$$

A series expansion with a computational package shows that all coefficients of the two sides  of (\ref{EtaRewriteTwo}) are equal up to the Sturm bound; hence the two sides are congruent modulo every prime $p$, i.e., the two sides are equal and the conjectured dissection holds. Theorem 1 follows immediately.  \hfill $\Box$

\subsection{Proofs of Theorem 2 and 3}

\noindent \textbf{Proof of Theorem 2:} Because $\frac{(q^{36};q^{36})_\infty}{(q^4;q^4)_\infty}$ is just the $q \rightarrow q^4$ blow-up of $\frac{(q^9;q^9)_\infty}{(q;q)_\infty}$, Theorem 2 holds if we additionally use the classical fact that 

$$\frac{(q^{12};q^{24})_\infty^2}{(q^4;q^8)_\infty^6} \equiv 1\, \text{ mod } \, 3$$

\noindent and that $b_9(16k-3)$ is even, which are Lemmas \ref{InsideThree} and \ref{EvenLemma}.

To prove Lemma \ref{EvenLemma}, we construct

$$C(z) := q^3 \left( \sum_{n=0}^\infty b_9(n) q^n \right) \cdot \left( \prod_{n=1}^\infty (1-q^n)^{64} \right) = \eta(9z) \eta^{63}(z).$$

We calculate that $C(z)$ is a modular form of weight 32 on $\Gamma_0(27)$.

Applying $T_2$ four times, we find that

$$(C \vert T_2^4)(z) \equiv q \left( \sum_{n=0}^\infty b_9(16n+13) q^n \right) \cdot \left( \prod_{n=1}^\infty (1-q^n)^4 \right) \, \text{ mod } \, 2.$$

The Sturm bound is $\frac{32}{12} \cdot 27 \cdot \frac{4}{3} = 96$, and computation verifies that all coefficients of $(C\vert T_2^4)(z)$ are even up to this bound.  Hence $(C\vert T_2^4)(z) \equiv 0$ mod 2.  Since $\prod_{n=1}^\infty (1-q^n)^4$ is invertible mod 2, we have that all coefficients of $\sum_{n=0}^\infty b_9(16n+13) q^n$ must be even.

Since both lemmas hold, $b_9(16n+13)$ must be divisible by 6, and Theorem 2 is proved. \hfill $\Box$

\phantom{.}

\noindent \textbf{Proof of Theorem 3:} The appearance of $B_9(q^4)$ in the proof of Theorem 2 is a surprisingly common feature in proofs of theorems for $\ell$-regular partitions and it is typical in papers on $\ell$-regular partitions to prove infinite families of congruences by means of such self-similarities.

In this case, by repeatedly applying the dissection to the 1 mod 4 subprogressions, and eventually taking the 3 mod 4 subprogression at some power of 4, we obtain Theorem 3 by induction. \hfill $\Box$

\section{Higher congruences}

The self-similarity phenomenon observed above appears for many $\ell$ and numerous small moduli.  It is easy to conjecture candidates for such self-similarity by numerical experimentation; so many exist that it seems an important problem in the area would be to develop a unified approach to the existence, frequency, formulas and symmetries among similarities of the form

$$\sum_{n=0}^{\infty} b_\ell(An+B) q^n \equiv q^j B_\ell(q^k) \, \text{ mod } m.$$

Cui and Gu, and Andrews, Hirschhorn and Sellers, prove such infinite families divisible by 2 (and a few by 3), using dissections; Furcy and Penniston prove infinite families modulo 3, using the properties of modular forms.  However, for $\ell$ and progression moduli of greater size, dissections begin to be unwieldy, and the modular form machinery almost exclusively gives progressions whose moduli are powers of the modulus of the congruence, or multiples thereof, which are only a small fraction of the congruences apparently holding.  

An alternative approach might be found in an interesting proof of a result by Penniston and Lovejoy (Theorem 3, \cite{FPenn}) from which it can be extracted that

\begin{corollary} $b_3(5n+2) \equiv b_3(7k+4) \equiv 0 $ mod 9 unless $n$ or $k$ are divisible by 5 or 7 respectively.\end{corollary}

Their proof is of an unusual sort for the literature, referencing the arithmetic geometry of a particular surface, and their remarks at the end of the paper pessimistic regarding future applications -- but a relation to the self-similarity problem suggests that there is some hope of use to gained by adapting that paper's arguments.  It seems very likely that this corollary can be considered a consequence of the self-similarity $\sum_{n=0} b_3(5n+2) q^n \equiv 2 B_3 (q^5)$ mod 9.

A congruence holding mod 9 brings us to a remaining area of related interest, to date apparently completely unexplored: the behavior of $l$-regular partitions modulo $2^j$ and $3^j$ for higher powers of $j$.  For instance, the $b_3(7n+4)$ progression mentioned above may hold mod 27.  There is also evidence for progressions in $b_9$ which are divisible by 27 and 54.  However, work of Gordon and Ono (\cite{GordonOno}) informs us that $b_9 (n)$ is almost always divisible by $3^j$ for any integer $j$, and so spurious progressions of this sort are possible, and caution advised.

Any congruence $p(An+B) \equiv 0$ mod $M$ for the ordinary partition function implies one mod $M$ for $b_{Ak}(An+B)$ for any multiple of $A$, and work of Ono and Mahlburg (with many others) gives us such congruences for powers of any prime at least 5.  However, we have no such congruences for the ordinary partition function with $A$ or $M = 2^j$ or $3^j$, and yet several such progressions exist for several $b_\ell$.  Therefore, it seems that an important open question is:

\phantom{.}

\noindent \textbf{Problem:} Do there exist progressions $b_\ell(An+B) \equiv 0$ mod $2^j$ or $3^j$ for indefinitely high $j$?

\end{document}